\documentclass[a4paper, 12pt]{amsart}

\usepackage[utf8]{inputenc}
\usepackage[T1]{fontenc}
\usepackage[english]{babel}

\usepackage{amssymb, amsmath, amsfonts, amscd}
\usepackage{amsthm}
\usepackage{bbm}

\usepackage[top=2cm, bottom=2cm, left=2cm, right=2cm]{geometry}

\newtheorem{lm}{Lemma}
\newtheorem{tm}{Theorem}
\newtheorem{cor}{Corollary}

\newcommand{\R}{\mathbb{R}}
\newcommand{\N}{\mathbb{N}}

\newcommand{\conv}{\mathrm{conv}}

\newcommand{\area}{\mathrm{area}\,}
\newcommand{\vol}{\mathrm{vol}\,}
\newcommand{\per}{\mathrm{per}\,}

\usepackage{hyperref}
\hypersetup{
    colorlinks=true,
    linkcolor=blue,
    citecolor=blue,
    urlcolor=blue
}

\begin{document}

\title[Inequalities for convex functions of random points]{Inequalities for convex functions of random points inside and on the boundary of convex bodies}

\author{A.~S.~Tokmachev}
\thanks{This work was supported by the Ministry of Science and Higher Education of the Russian Federation (agreement 075-15-2025-344 dated 29/04/2025 for Saint Petersburg Leonhard Euler International Mathematical Institute at PDMI RAS)}
\thanks{The work was supported by the Theoretical Physics and Mathematics Advancement Foundation ``BASIS''}
\thanks{Supported in part by the Moebius Contest Foundation for Young Scientists}
\address{St.~Petersburg Department of Steklov Institute of Mathematics, St.~Petersburg, Russia}
\email{chief.tokma4eff@yandex.ru}


\begin{abstract}
    Let $K$ be a convex body in $\R^d$, and let $I$ and $B$ be random points uniformly distributed inside $K$ and on its boundary, respectively. We prove that if $d=2$ and $\mathbb E I = \mathbb E B$, or if $K$ is a circumscribed polytope with the center of the inscribed sphere coinciding with $\mathbb E I = \mathbb E B$, then $I$ is dominated by $B$ in the convex order. As a consequence, for any function $\varphi$ convex in each argument, the expectation $\mathbb E \varphi(I_1,\dots,I_k)$ does not exceed $\mathbb E \varphi(B_1,\dots,B_k)$. This yields, in particular, an inequality between the moments of random chords $\mathbb E |I_1 - I_2|^p \leqslant \mathbb E |B_1 - B_2|^p$ for all $p \geqslant 1$, confirming the Zaporozhets--Tarasov conjecture for the indicated class of bodies, and extends to inequalities for mean volumes of random simplices.
\end{abstract}

\maketitle

\section{Introduction}

Let $K$ be a convex figure in the plane. By this we mean that $K$ is a convex compact set with nonempty interior; in the sequel we refer to such sets as convex bodies, so as not to introduce separate terminology for the higher-dimensional case.

It is convenient to fix notation for the random points used throughout the paper. Let $I$ be uniformly distributed inside $K$ (with respect to the Lebesgue measure), and let $B$ be uniformly distributed on the boundary $\partial K$ (with respect to the surface measure). Independent copies of these variables will be denoted by $I_1, I_2, \dots$ and $B_1, B_2, \dots$, respectively.

\vskip 6pt

In 1864, Sylvester posed the following problem \cite{Sylvestr}: four points $I_1, I_2, I_3, I_4$ are chosen uniformly and independently inside $K$. What is the probability that their convex hull $\conv(I_1, I_2, I_3, I_4)$ is a triangle? Clearly, this probability depends on the shape of $K$. In 1918, Blaschke \cite{Blaschke} proved that for any planar convex body $K$,
\begin{equation*}
    \frac{35}{12 \pi^2} \leqslant \mathbb P(\conv(I_1, I_2, I_3, I_4) \text{ is a triangle}) \leqslant \frac{1}{3},
\end{equation*}
where the lower bound is attained on ellipses and the upper bound on triangles.

This problem admits a reformulation in terms of areas. Let $\area(\cdot)$ denote the area. A simple computation shows that
\begin{align*}
    \mathbb P(\conv(I_1, I_2, I_3, I_4) \text{ is a triangle}) &= 4 \mathbb P(I_4 \in \conv(I_1, I_2, I_3)) = 4\frac{\mathbb E\, \area(\conv(I_1 I_2 I_3))}{\area K}.
\end{align*}
Set $\alpha(K) = \mathbb E\, \area(\conv(I_1 I_2 I_3))$. Blaschke's inequality then becomes
\begin{align*}
   \frac{35}{48\pi^2}\leq\frac{\alpha(K)}{\area K}\leq\frac{1}{12},
\end{align*}
which gives sharp bounds for the normalized mean area of a random triangle with vertices uniformly distributed inside $K$.

\vskip 6pt

In \cite{BGTZ}, the analogous problem for two points was considered: the length of a random segment with endpoints chosen uniformly and independently inside $K$, normalized by the perimeter of $K$, was estimated. It was shown that for
\begin{align*}
   \Delta(K) = \mathbb{E}\,|I_1 - I_2|
\end{align*}
the following sharp inequalities hold:
\begin{align}\label{0020}
   \frac{7}{60}<\frac{\Delta(K) }{\per K}<\frac{1}{6}.
\end{align}
This result was also extended to higher dimensions, where $\Delta$ is normalized by the mean width of the body. In dimension~2, by Cauchy's formula (see, e.g., \cite{Cauchy}), the mean width coincides with the perimeter multiplied by $1/\pi$.

Although both bounds in~\eqref{0020} are strict, they are optimal: the lower bound is asymptotically attained on a sequence of isosceles triangles whose base is fixed while the height tends to zero, and the upper bound is asymptotically attained on a sequence of rectangles degenerating to a segment (which itself has empty interior and does not belong to the class of convex bodies considered here).

\vskip 6pt

Now consider the analogous functionals with points taken on the \emph{boundary} of the body:
\begin{align}\label{1024}
    \theta(K) &= \mathbb{E}\,|B_1 - B_2|,\\
    \psi(K) &= \mathbb{E}\,\area(\conv(B_1, B_2, B_3)).
\end{align}

\vskip 6pt

Zaporozhets and Tarasov conjectured that for every convex body in $\R^d$,
\begin{align}\label{0156}
\Delta(K)<\theta(K).
\end{align}
To the best of our knowledge, this conjecture remains open even in dimension $d=2$.

Gusakova and Zaporozhets proposed the following possible approach. They suggested that for any planar convex body,
\begin{align}\label{0024}
   \frac{1}{6}<\frac{\theta(K) }{\per K}\leq \frac{2}{\pi^2},
\end{align}
where equality in the upper bound holds if and only if $K$ is a disk, while the lower bound is unattainable in the class of convex bodies but is reached on a segment, which can be approximated by rectangles. Clearly, the lower bound in~\eqref{0024} together with the upper bound in~\eqref{0020} would immediately imply~\eqref{0156}.

\vskip 6pt

In the author's works \cite{Tokmachev1, Tokmachev2}, the upper bounds in conjecture~\eqref{0024} and related results for $\psi$ were proved. Namely, the disk maximizes both the normalized mean distance between two boundary points and the normalized mean area of a triangle with three boundary vertices, for a fixed perimeter.

\vskip 6pt

In 2025, A.~S.~Lotnikov \cite{Lotnikov} obtained a number of results related to conjecture~\eqref{0156} for centrally symmetric planar bodies. He showed that for such bodies not only does $\Delta(K) < \theta(K)$ hold, but a stronger stochastic domination statement is true: for any $n>1$, $\mathbb{E}|B_1 B_2|^n \geqslant \mathbb{E}|I_1 I_2|^n$. He also established exact formulas relating moments of interior and boundary chords for circumscribed polygons.

\medskip

\textbf{The aim of the present paper} is to extend the chord-length inequality to a wider class of bodies. We show that in the plane, the inequality $\mathbb{E}|I_1 I_2|^p \leqslant \mathbb{E}|B_1 B_2|^p$ for all $p \geqslant 1$ holds under the single condition $\mathbb E I = \mathbb E B$ (which is automatically true for centrally symmetric bodies). Moreover, we prove a substantially more general fact: the inequality holds for arbitrary functions convex in each argument. As a corollary, we obtain an analogous inequality for the mean volume of a random simplex with vertices chosen inside and on the boundary of the body.

\section{Main results}

Our main result is the following theorem.

\begin{tm}\label{main_tm}
    Let $K$ be a convex body in $\R^d$ satisfying one of the following:
    \begin{itemize}
        \item $d = 2$ and $\mathbb E I = \mathbb E B$;
        \item $K$ is a circumscribed polytope whose inscribed sphere center $O$ coincides with $\mathbb EI$ and $\mathbb EB$.
    \end{itemize}
    Let $\varphi: (\R^d)^k \to \R$ be a function convex in each argument.
    Then
    \begin{align}
        \mathbb E \varphi(I_1, \ldots, I_k) \leqslant  \mathbb E \varphi(B_1, \ldots, B_k),
    \end{align}
    where $I_1,\dots,I_k$ and $B_1,\dots,B_k$ are independent samples uniformly distributed inside $K$ and on $\partial K$, respectively.
\end{tm}

\vskip 6pt

Taking $\varphi(x, y) = |x - y|^p$ ($p \geqslant 1$) gives a generalization of Lotnikov's result.

\begin{cor}\label{cor:distance}
    Under the conditions of Theorem~\ref{main_tm}, for any $p \geqslant 1$,
    \begin{align}
        \mathbb E |I_1 - I_2|^p \leqslant \mathbb E|B_1 - B_2|^p.
    \end{align}
\end{cor}

For $p = 1$, Corollary~\ref{cor:distance} confirms the Zaporozhets--Tarasov conjecture for the bodies in question.

\vskip 6pt

Another example of a function convex in each argument is
\[
\varphi(x_1,\ldots, x_n) = (\vol_d(\conv\{x_1,\ldots, x_n\}))^p,
\qquad p \geqslant 1.
\]

\begin{cor}\label{cor:volume}
    Under the conditions of Theorem~\ref{main_tm}, for any $n \in \N$ and $p \geqslant 1$,
    \begin{align}
        \mathbb E \bigl[\vol_d(\conv\{I_1,\ldots, I_n\})^p\bigr] 
        \leqslant 
        \mathbb E\bigl[\vol_d(\conv\{B_1,\ldots, B_n\})^p\bigr].
    \end{align}
\end{cor}

For $n = d+1$ and $p = 1$, this yields $\alpha(K) \leqslant \psi(K)$ for the mean volumes of random simplices with vertices inside and on the boundary.

\section{Reduction to the convex order}

The key ingredient in the proof of Theorem~\ref{main_tm} is the convex order between the distributions of $I$ and $B$. Recall that $X \preceq_{cx} Y$ (in the convex order) if $\mathbb E \varphi(X) \leqslant \mathbb E \varphi(Y)$ for every convex function $\varphi$. By Strassen's theorem \cite{Strassen}, $X \preceq_{cx} Y$ is equivalent to the existence of a coupling with $\mathbb E[Y \mid X] = X$ almost surely; in particular, $\mathbb E X = \mathbb E Y$ is necessary.

\vskip 6pt

To prove Theorem~\ref{main_tm} it suffices to establish $I \preceq_{cx} B$; the full statement then follows from independence and properties of conditional expectations. Thus we need only prove:

\begin{tm}[Convex order for $I$ and $B$]\label{main_tm_new}
    Under the same assumptions on $K$ as in Theorem~\ref{main_tm},
    \begin{align}
        I \preceq_{cx} B,
    \end{align}
    i.e., $\mathbb E \varphi(I) \leqslant \mathbb E \varphi(B)$ for every convex $\varphi \colon \R^d \to \R$.
\end{tm}

\vskip 6pt

The case of a circumscribed polytope is due to Pasteczka \cite{Pastechka} (Theorem~2). The next section contains the proof for $d=2$, which is the main contribution of this note.

\section{Proof of Theorem~\ref{main_tm_new} for \texorpdfstring{$d=2$}{d=2}}

We need two auxiliary facts: the Ohlin lemma and a characterization of the convex order via one-dimensional projections.

\begin{lm}[Ohlin \cite{Ohlin}]
    Let $X$ and $Y$ be integrable random variables with $\mathbb E X = \mathbb E Y$. If their distribution functions $F_X$ and $F_Y$ cross exactly once, with $F_X \leqslant F_Y$ to the left of the crossing point and $F_X \geqslant F_Y$ to the right, then $\mathbb E f(X) \leqslant \mathbb E f(Y)$ for every convex $f \colon \R \to \R$.
\end{lm}

\vskip 6pt

The next fact is well known (see, e.g., \cite{MullerStoyan}).

\begin{lm}\label{stat_conv}
    For $X, Y \in \R^d$, one has $X \preceq_{cx} Y$ if and only if for every unit vector $u \in S^{d-1}$ and every convex $f \colon \R \to \R$,
    \begin{align}
        \mathbb E f(\langle X, u\rangle) \leqslant \mathbb E f(\langle Y, u\rangle),
    \end{align}
    where $\langle \cdot, \cdot \rangle$ is the scalar product.
\end{lm}

\medskip

Now we prove the theorem. Let $d = 2$ and $\mathbb E I = \mathbb E B$. Fix a unit vector $u \in S^1$ and let $I_u = \langle I, u\rangle$, $B_u = \langle B, u\rangle$ be the projections onto the line spanned by $u$. We verify the Ohlin condition for $I_u$ and $B_u$.

First, $\mathbb E I = \mathbb E B$ implies
\begin{align*}
    \mathbb E I_u = \langle \mathbb E I, u\rangle = \langle \mathbb E B, u\rangle = \mathbb E B_u.
\end{align*}

Second, by Lemma~2 of \cite{BartheEtAl}, the distribution functions $F_{I_u}$ and $F_{B_u}$ cross exactly once, with the sign pattern required in Ohlin's lemma (the crossing point $x_0$ is the projection of the ``widest part'' of $K$ onto this line). Hence,
\begin{align*}
    \mathbb E f(I_u) \leqslant \mathbb E f(B_u)
\end{align*}
for every convex $f$. By Lemma~\ref{stat_conv}, this yields $I \preceq_{cx} B$, completing the proof of the two-dimensional case.


\begin{thebibliography}{99}

\bibitem{Sylvestr} J.~J.~Sylvester. Problem 1491. \emph{The Educational Times}, April 1864.

\bibitem{Blaschke} W.~Blaschke. \"Uber affine Geometrie XI: L\"osung des ``Vierpunktproblems'' von Sylvester aus der Theorie der geometrischen Wahrscheinlichkeiten. \emph{Ber. Verh. S\"achs. Akad. Wiss. Leipzig, Math.-Phys. Kl.}, 69:436--453, 1917.

\bibitem{BGTZ} G.~Bonnet, A.~Gusakova, Ch.~Th\"ale, and D.~Zaporozhets. Sharp inequalities for the mean distance of random points in convex bodies. \emph{Adv. Math.}, 386:107799, 2021.

\bibitem{Cauchy} R.~V.~Ambartzumian. \emph{Combinatorial Integral Geometry with Applications to Mathematical Stereology}. John Wiley \& Sons, Chichester, 1982.

\bibitem{Tokmachev1} A.~S.~Tokmachev. Mean distance between random points on the boundary of a convex body. \emph{J. Math. Sci. (N. Y.)}, 286(5):798--806, 2024. [Translation from: Zap. Nauchn. Sem. POMI, 510:248--261, 2022.]

\bibitem{Tokmachev2} A.~S.~Tokmachev. On the mean area of a triangle inscribed in a convex figure. \emph{Zap. Nauchn. Sem. POMI}, 525:134--149, 2023. (In Russian.)

\bibitem{Lotnikov} A.~S.~Lotnikov. Mean distance between random points in centrally symmetric convex bodies. \emph{Zap. Nauchn. Sem. POMI}, 544:211--235, 2025. (In Russian.)

\bibitem{Strassen} V.~Strassen. The existence of probability measures with given marginals. \emph{Ann. Math. Statist.}, 36(2):423--439, 1965.

\bibitem{MullerStoyan} A.~M\"uller and D.~Stoyan. \emph{Comparison methods for stochastic models and risks}. Wiley Series in Probability and Statistics. John Wiley \& Sons, Chichester, 2002.

\bibitem{Pastechka} P.~Pasteczka. Jensen-type geometric shapes. \emph{Ann. Univ. Paedagog. Crac. Stud. Math.}, 19:103--111, 2020.

\bibitem{Ohlin} J.~Ohlin. On a class of measures of dispersion with application to optimal reinsurance. \emph{Astin Bulletin}, 5(2):249--266, 1969.

\bibitem{BartheEtAl} F.~Nazarov, D.~Ryabogin, and V.~Yaskin. On the maximal distance between the centers of mass of a planar convex body and its boundary. \emph{Discrete Comput. Geom.}, 73:1016--1036, 2025.

\end{thebibliography}
\end{document}